\documentclass[11pt]{article}
\usepackage{amsmath,amssymb,amsthm,xcolor,graphicx,bbm}
\usepackage{accents}
\usepackage[breaklinks=true, colorlinks=true, urlcolor=cyan]{hyperref}
\usepackage{esint}
\usepackage{listings}
\usepackage{enumitem}

\usepackage[top=1in, bottom=1in, left=1.25in, right=1.25in, marginparwidth=1.1in, marginparsep=0.05in]{geometry}

\usepackage{titletoc}
\titlecontents{section}[1.5em]{\vspace{-2pt}}%
{\thecontentslabel\enspace\enspace}%
{}%
{\titlerule*[0.5pc]{.}\contentspage}
\titlecontents{subsection}[3em]{\vspace{-2.5pt}}%
{\thecontentslabel\enspace\enspace}%
{}%
{\titlerule*[0.5pc]{.}\contentspage}%

\numberwithin{equation}{section}
\newtheorem{theorem}{Theorem}[section]

\newtheorem{lemma}[theorem]{Lemma}
\newtheorem{proposition}[theorem]{Proposition}

\theoremstyle{remark}
\newtheorem{remark}[theorem]{Remark}
\newcommand{\bke}[1]{\left ( #1 \right )}
\newcommand{\bkt}[1]{\left [ #1 \right ]}
\newcommand{\bket}[1]{\left \{ #1 \right \}}
\newcommand{\norm}[1]{\left  \| #1  \right \|}
\newcommand{\bka}[1]{{\langle #1 \rangle}}

\newcommand\al{\alpha}

\newcommand\ga{\gamma}
\newcommand\de{\delta}
\newcommand\ep{\epsilon}

\newcommand\e {\varepsilon}

\renewcommand\th{\theta}

\newcommand\si{\sigma}
\newcommand\om{\omega}

\newcommand\De{\Delta}

\newcommand\Om{\Omega}
\newcommand{\R}{\mathbb{R}}
\newcommand{\RR}{\mathbb{R}}

\newcommand{\FF}{\mathbb{F}}
\newcommand{\GG}{\mathbb{G}}

\newcommand{\NN}{\mathbb{N}}

\newcommand{\cW}{\mathcal{W}}

\renewcommand{\div}{\mathop{\rm div}\nolimits}

\newcommand{\supp} {\mathop{\mathrm{supp}}\nolimits}

\newcommand{\dist} {\mathop{\mathrm{dist}}\nolimits}

\newcommand{\pd}{\partial}
\newcommand{\nb}{\nabla}

\newcommand{\lec}{{\ \lesssim \ }}

\newcommand{\oo}{\infty}

\newcommand{\EQ}[1]{\begin{equation} #1 \end{equation}}

\newcommand{\EQS}[1]{\begin{equation}\begin{split} #1 \end{split}\end{equation}}
\newcommand{\EQN}[1]{\begin{equation*}\begin{split} #1 \end{split}\end{equation*}}
\newcommand{\EN}[1]{\begin{enumerate}[label={\textup{(\roman{enumi})}}] #1 \end{enumerate}}

\usepackage{mathtools} %

\newcommand{\loc}{\mathrm{loc}}

\newcommand{\arxiv}[1]{\href{http://arxiv.org/abs/#1}{arXiv:#1}}

\begin{document}
\title{Existence and regularity for perturbed Stokes system with critical drift in 2D}
\author{Misha Chernobai and Tai-Peng Tsai}
\date{}
\maketitle
\begin{abstract}
We consider a perturbed Stokes system with critical divergence-free drift in a bounded Lipschitz domain in $\R^2$, with sufficiently small Lipschitz constant $L$.
It extends our previous work in $\R^n$, $n\ge 3$, to two-dimensional case.
For large drift in weak $L^2$ space, we prove unique existence of $q$-weak solutions for forces in $L^q$ with $q$ close to $2$.
The dimension 2 case requires special care as the energy estimate for $2$-weak solutions is not defined.
We first prove the existence of $2$-weak solutions by approximations without uniqueness. We then show the results for $q$ slightly larger than 2. We finally get the results for $q<2$ by duality, which implies uniqueness for $q=2$.
We also show for any $1<q<\infty$ the unique existence of $q$-weak solutions for forces in $L^q$ when the drift is sufficiently small in weak $L^2$. Using similar methods we can also prove analogous results for scalar equations with
divergence-free drifts in weak $L^{2}$ space.

\medskip

\emph{2020 Mathematics Subject Classifications}: 35J47, 35Q35
\end{abstract}

\setcounter{tocdepth}{1}
\tableofcontents

\section{Introduction}
\label{sec1}

For a bounded Lipschitz domain $\Om$  in $\R^2$ we consider the solutions of the perturbed Stokes system $(u,\pi):\Om \to \R^2 \times \R$ that satisfy
\EQ{\label{drifteq}
-\Delta u +b\cdot\nabla u+\nabla \pi=\div \GG,\quad \div u=0, \quad u|_{\pd \Om}=0,
}
with given drift $b$ in the critical space, $b\in L^{2,\oo}(\Omega)^2$, $\div b=0$, and  force $\GG\in L^q(\Om)^{2\times 2}$, $1 <q<\oo$. We denote by
$L^{q,r}(\Om)$ the Lorentz spaces. For derivatives of vectors and matrices we use the following notations
$(\div \GG)_i=\pd_j \GG_{ji}$ and $(\nb \zeta)_{ji} = \pd_j \zeta_i$. Since $\div b=0$  we can integrate by parts in the following terms
\EQ{\label{1.2}
\De u = \div (\nb u),\quad
b\cdot\nabla u = \div (b\otimes u), \quad \int_\Om  \div \GG \cdot \zeta = -
\int_\Om  \GG:\nb \zeta.
}
For zero divergence spaces or spaces with zero average we use the following notations
(for $\Om\subset\R^n$, $n\ge2$)
\EQS{
C^{\infty}_{c}(\Omega) &= \big\{ u \in C^{\infty}(\Omega):~ \supp\{u\}\Subset \Omega\big\},\\
C^{\infty}_{c,\si}(\Omega) &= \big\{ u \in C^{\infty}_c(\Omega)^n:\  \div u=0\big\},\\
W^{1,q}_{0,\si}(\Omega) &= \big\{ u \in W^{1,q}_0(\Omega)^n:\  \div u=0\big\},\\ \quad L^q_0(\Om)&=\bket{\pi \in L^q(\Om):\  \textstyle \int_\Om \pi \,dx=0}.
}

Next we consider $(u,\pi)$ as a \emph{weak solution pair} of \eqref{drifteq} if $u\in W^{1,q}_{0,\si}(\Omega)$, $\pi\in L^q_0(\Omega)$, and $(u,\pi)$ satisfies  \eqref{drifteq} in distributional sense,\EQ{\label{soln-pair}
\int_{\Om} (\nb u - b \otimes u): \nb \zeta - \pi \div \zeta  = -\int_\Om \GG : \nb \zeta,\quad
\forall \zeta\in C_c^{\infty}(\Om)^2.
}
Here we used \eqref{1.2} to integrate by parts.
We say $u$ is a \emph{$q$-weak solution} of \eqref{drifteq} if $u\in W^{1,q}_{0,\si}(\Omega)$ and $u$ satisfies \eqref{drifteq} in weak sense, i.e.,
\EQ{\label{weak-soln}
\int_{\Om} (\nb u - b \otimes u): \nb \zeta   = -\int_\Om \GG : \nb \zeta,\quad
\forall \zeta\in C_{c,\sigma}^{\infty}(\Omega).
}
In case $q=2$ we call $u$ a \emph{weak solution}. In this paper we also consider a dual system of \eqref{drifteq}
\EQ{\label{drifteq_dual}
-\Delta v -\nabla\cdot(b\otimes v)+\nabla \pi=\div \FF,\quad \div v=0, \quad v|_{\pd \Om}=0.
}
This system is equivalent to \eqref{drifteq} when the drift is divergence free, $\div b=0$, but in general case it is different.

The boundary value problem of the perturbed Stokes system \eqref{drifteq} and the corresponding scalar equation \eqref{scalar.eq}
has a huge amount of literature when we consider subcritical drift $b\in L^p(\Om)$, $\Om \subset \R^n$, with $p> n \ge2$. In the subcritical case we may treat the drift term $b\cdot \nb u$ as a perturbation to $\De u$.

The critical case is when $b\in L^n(\Om)$ or $b\in L^{n,\infty}(\Om)$. System \eqref{drifteq}
with critical drift $b\in L^n(\Om)$, $n\ge2$, has previously been studied;
See Dindoš and Mitrea \cite[Theorem 6.1]{DinMit}, Kim \cite[\S4]{Kim09}, Choe \& Kim \cite[Theorem 18]{ChoeKim}, and
Amrouche \& Rodríguez-Bellido \cite[\S5]{AmRB}.
They contain estimates similar to \eqref{main1-eq}.
For general results on scalar equation \eqref{scalar.eq} with drift  $b \in L^n$, see
Kim \& Kim \cite{KK}.

When one considers $b\in L^{n,\infty}(\Om)$, the case $n \ge 3$ is considered by the authors
in \cite[Theorems 1.4, 1.5]{CT1}. They correspond to Proposition \ref{prop2}
and Theorem \ref{main1} below for dimension $n=2$.

In this paper we extend our previous work for higher dimensions $n \ge 3$ \cite{CT1} to $n=2$. Our main theorem below extends \cite[Theorem 1.5]{CT1} to the dimension 2 case and proves existence of the solutions in $W^{1,q}$ for arbitrary  drift in weak $L^2$ with $q$ close to $2$.
\begin{theorem}\label{main1}
Let $\Om$ be a bounded Lipschitz domain in $\R^2$ with Lipschitz constant $L < L_0$, where $0<L_0<1$ is a sufficiently small constant. Assume $b\in L^{2,\infty}(\Omega)^2$ and $\div b=0$.
Then there exists $p_0(\Om, \|b\|_{L^{2,\infty}})>2$ such that for any $p_0'<q<p_0$ and $\GG\in L^q(\Om)^{2\times 2}$, there exists a  weak solution pair $u\in W^{1,q}_{0,\si}(\Omega)$ and $\pi\in L^q_0(\Omega)$ of \eqref{drifteq} such that
\EQ{\label{main1-eq}
\|u\|_{ W^{1,q}(\Omega)}+\|\pi\|_{L^q(\Omega)}\le C(\Omega,\|b\|_{L^{2,\infty}})\|\GG\|_{L^q(\Omega)}.
}
For $p_0'<q<p_0$, the above pair $(u,\pi)$ is the unique weak solution pair  in $W^{1,q}_0\times L^q_0(\Om)$ without assuming \eqref{main1-eq}.
\end{theorem}

\emph{Comments on Theorem \ref{main1}:}
\EN{
\item This theorem extends the higher dimensional result \cite[Theorem 1.5]{CT1} to the dimension 2 case.
In \cite[Theorem 1.5]{CT1} we take $L_0=1/2$. In this paper $L_0$ is specified after \eqref{change-var} and after \eqref{pressure_bound_est}, and is computable.

\item  Since $\Om$ is Lipschitz we can assume that there exists $R_0>0$ such that for any $x_0\in \partial\Om$ we have
\EQ{
\Om\cap B_{2R_0}(x_0)=\{(x',x_n)\in B_{2R_0}(0),x_n>\gamma(x')\}.
}
The constant $C$ in \eqref{main1-eq} will actually depend only on $R_0$, $\text{diam}~\Om$, and $\|b\|_{L^{2,\infty}}$.

\item If we further assume $b \in L^2(\Om)$ for dimension $n=2$,
the unique existence of weak solution in $W^{1,q}$ of  \eqref{drifteq}, $1<q<\infty$, was shown in
Dindoš and Mitrea \cite[Theorem 6.1]{DinMit}.

\item We may simply assume $\Om$ is a bounded $C^1$ domain, as for
a bounded $C^1$ domain and any constant $\lambda>0$, each point $x$ on $\pd\Om$ has a neighborhood $U_x$ such that $\pd \Om \cap U_x$ is the graph of a Lipschitz function with a Lipschitz constant less than $\lambda$.

}

We also formulate the following result for any $1<q<\infty$ under a smallness assumption on $\|b\|_{L^{2,\infty}}$.

\begin{proposition}\label{prop2}
  Let $1<q<\infty$ and $\Om\subset\RR^2$ be a bounded Lipschitz domain with sufficiently small  ($q$-dependent) Lipschitz constant $L>0$. There exists $\varepsilon=\e(q,\Om)>0$ such that for any $b\in L^{2,\infty}(\Om)^2$ with $\div b=0,~\|b\|_{L^{2,\infty}}\le\varepsilon$ and any $\GG \in L^q(\Om)^{2\times 2}$ there exists a unique $q$-weak solution  pair $u\in W^{1,q}_{0,\si}(\Omega)$ and $\pi\in L^q(\Omega)$ of \eqref{drifteq} such that $\int_\Om \pi\,dx=0$ and
\EQ{
\|u\|_{ W^{1,q}(\Omega)}+\|\pi\|_{L^q(\Omega)}\le C\|\Bbb G\|_{L^q(\Omega)}.
}
\end{proposition}

In the proposition, the case $q=2$ is already included in Theorem \ref{main1}. For $q\not =2$,
the proof is the same as \cite[Theorem 1.4]{CT1}: The a priori bound of $b\cdot \nb u$ in $W^{-1,q}(\Om)$ of  \cite[Lemma 2.4]{CT1} is valid for $n=2 \not = q$, and one can use the method of continuity (and duality for $q<2$) to construct the solution when $\|b\|_{L^{2,\infty}}$ is sufficiently small.
We skip the details.

For corresponding results on scalar equations
\EQ{\label{scalar.eq}
-\Delta u +b\cdot\nabla u=\div  G, \quad u|_{\pd \Om}=0,
}
with critical drift coefficient $b \in L^n$, see
Kim \& Kim \cite{KK}. When $b \in L^{n,\infty}(\Om)$ for $n \ge 3$, see \cite{KiTs,KPT} which also cover more general assumptions on drift $b$ as well as %
\cite{CS20,CS24} and Kwon \cite{Kwon} for $n=2$ case or drift in Morrey spaces which include $L^{n,\infty}$.
There are also results on H\"older continuity of the solutions, for example with $b\in BMO^{-1}(\R^3)$, $\div b=0$, see Seregin-Sverak-Silvestre-Zlatos \cite{SSSZ}.

Using the same proof scheme of Theorem \ref{main1} we can also prove the following theorem for the  scalar equation \eqref{scalar.eq}.
\begin{theorem}\label{th1.2}
  Let $\Om$ be a bounded Lipschtiz domain in $\RR^2$.  Assume $b\in L^{2,\infty}(\Omega)$ and $\div b=0$.
Then there exists $p_0(\Om,\|b\|_{L^{2,\infty}})>2$ such that for any $p_0'<q<p_0$ and $G\in L^q(\Om)$, there exists a  weak solution pair $u\in W^{1,q}_{0,\si}(\Omega)$ of \eqref{scalar.eq} such that
\EQ{\label{main2-eq}
\|u\|_{ W^{1,q}(\Omega)}\le C(\Omega,\|b\|_{L^{2,\infty}})\|G\|_{L^q(\Omega)}.
}
For $p_0'<q<p_0$, the above solution $u$ is the unique weak solution in $W^{1,q}_0(\Om)$ without assuming \eqref{main2-eq}.
\end{theorem}

This theorem improves $\alpha=0$ case of \cite[Theorem 1.1]{CS20} and $n=2$ case of Kwon \cite[Theorem 1.1]{Kwon} which show that for any $q>2$, there is $p \in (2,q)$ depending on $q$, $\norm{b}_{L^{2,\infty}}$ and $\Om$ such that for any $G \in L^q(\Om)$, there is a unique $p$-weak solution $u$ of \eqref{scalar.eq} with force $G$. Moreover,
\[
\norm{\nb u}_{L^p(\Om)} \le C \norm{G}_{L^q(\Om)}.
\]
Theorem \ref{th1.2} improves this result by asserting that $p=q$ when $2<q<p_0$ is sufficiently close to $2$. This improvement further allows us to prove a priori bounds for $q$-weak solutions for $q \in (p_0',2]$ by duality, which yields existence and uniqueness
 for $q \in (p_0',2]$.

We skip the proof of Theorem \ref{th1.2} since it is similar to Theorem \ref{main1} and is easier, without the need of pressure estimate. Because it has no pressure estimate, the Lipschitz constant of the boundary $\pd \Om$ is allowed to be arbitrarily large.

\medskip

We now explain the new difficulties and ideas.
Theorem \ref{main1}
extends the higher dimensional result \cite[Theorem 1.5]{CT1} to the 2D case, and
we will still use Gehring's approach. In this approach, we need to do local energy estimate in every ball and prove a reversed H\"older inequality in a general ball with a uniform constant.

The first new difficulty is that we cannot use $\zeta = u \eta^4$ for some cut-off function $\eta$ as the test function in the weak form \eqref{soln-pair}, because $b \otimes u$ is no longer in $L^2$ as in the 3D case, and
 the term
\EQ{\label{drift-integral}
\int_\Om b \otimes u : \nb( u \eta^4)
}
is not integrable. Hence we will mollify the drift $b$, so that we can derive an a priori bound. Since there is no theorem that extends a div-free vector field in $L^{2,\infty}(\Omega)$ to $L^{2,\infty}(\R^2)$ for a general bounded domain, we will consider the approximation problems in subdomains $\Om_n$ that converge to $\Om$.

The second new difficulty comes from the lower integrability of the pressure.
We will
bound the pressure using Wolf's local pressure projection with a uniform constant as in \cite{CT1}.
Denote $\Omega_R=\Omega\cap B_R(x_0)$ for some $x_0\in \overline \Om$. In \cite{CT1}, it suffices to bound $\pi -(\pi)_{\Omega_R}$ in $L^2(\Omega_R)$. However,  for dimension 2  we need to bound $\pi -(\pi)_{\Omega\cap B_R}$ in $L^q(\Omega_R)$, $q<2$. (We choose $q=4/3$ for simplicity.)
We have to choose $q<2$ because
\[
\pi \lec \nb u + b \otimes u -\GG
\]
and $b \otimes u$ is at best in $L^{2-}$ for $u \in W^{1,2}_\loc$ even if we assume $b \in L^2$. This would not be an issue if $b \in L^{r}$, $r>2$, and we can only assume $u \in W^{1,2}_\loc$ which is the assumption of the reversed H\"older inequality.

To obtain a uniform constant of Wolf's local pressure projection in $L^q(\Omega_R)$ for $q<2$ gives a significant higher restriction on the regularity of the domain. Instead of a Lipschitz domain with Lipschitz constant $L$ less than 1/2, we have to assume $L$ is sufficiently small and do boundary stretching in the proof for the 2D case.

The third difficulty is the approximation of a Lipschitz domain. In order to construct the solution we are using an approximation method, therefore the first logical idea would be to construct a sequence of $b_n\in L^{2,\infty}(\Om)$, $\div b=0$, such that $b_n$ converge in a sense that allows us to prove that the limit will satisfy the equation, for example convergence almost everywhere is sufficient and $\|b_n\|_{L^{2,\infty}(\Om)}$ are bounded uniformly. For example in star-shaped domains it is trivial, but for arbitrary Lipschitz domains this question is not so clear. To deal with this issue, following Verchota \cite{Verchota_1,Verchota_2} and later Kwon \cite{Kwon}, Chernobai-Shilkin \cite{CS24}, we instead approximate the domain by Lipschtiz sub-domains with controlled Lipschitz constant and on sub-domains we can use standard mollification of $b$ and construct an approximation sequence as zero extension of solutions in these sub-domains. See Section \ref{sec6} for details.

The rest of the paper is organized as the following. In section 2 we have preliminary results concerning Lorentz spaces and the pressure, which will be used in the later sections. In section 3 we derive the existence of weak solutions in $W^{1,2}(\Omega)$. Section 4 contains Wolf's local pressure projection method that will be used to estimate the pressure. In section 5, assuming higher integrability of the drift, we prove uniform local estimate for the gradient which will imply higher integrability due to Gehring's lemma. Section 6 contains proof of the main Theorem \ref{main1}. In the Appendix we put a proof of geometrical lemma that was used in section 4.

\section{Preliminaries}\label{sec2}
In this section we give a few preliminary results. We denote the average of a function $f$ over a bounded open set $E\subset \R^n$ as
\[
(f)_E  = \frac 1{|E|} \int_E f\,dx = \fint_E f\,dx.
\]
For an open set $\Om\subset \R^n$ we denote
\[
\Om_r(x_0)=\Om_{x_0,r} = \Om \cap B_r(x_0).
\]
We use the standard notations for the conjugate Sobolev space $W^{1,-q}(\Om)^m$, $1<q<\infty$, $m\in\NN$, for the bounded $\Om$  with the following norm
\EQ{
\norm{f}_{W^{1,-q}(\Om)^m} = \sup _{u \in W^{1,q'}_0(\Om)^m,\,\norm{\nb u}_{L^{q'}}(\Om)\le1}
\bka{f, u}.
}

For our energy estimates we will need H\"older and
Sobolev inequalities in Lorentz spaces. R. O’Neil \cite{ON} proved H\"older inequality in Lorentz spaces, and for cases $p \le 1$ or $q< 1$ see \cite[Lemma 3.1]{KPT} .

\begin{lemma}[H\"older inequality in Lorentz spaces]\label{Holder}
Let $\Om$ be any domain in $\R^n$, then for any $0< p,p_1,p_2< \infty$ and $0< q,q_1,q_2\le\infty$ that satisfy
$$
\frac{1}{p}=\frac{1}{p_1}+\frac{1}{p_2}\quad\text{and}\quad \frac{1}{q} \le  \frac{1}{q_1}+\frac{1}{q_2},
$$
there is a constant $C=C(  p_1 , p_2 ,   q_1 , q_2 ,q )>0$ such that
$$
\|fg\|_{p,q} \le C  \|f\|_{L^{p_1,q_1}(\Om)} \|g\|_{L^{p_2,q_2}(\Om)}
$$
for all $f \in L^{p_1,q_1}(\Omega ) $ and $g \in L^{p_2,q_2}(\Omega )$.
\end{lemma}

Next we formulate Sobolev inequality in Lorentz spaces. The proof follows from
Sobolev inequality in Lorentz spaces in the whole $\R^n$ (see \cite[Remark 7.29]{AdamsFournier} and \cite{PO}), and the extension theorem from $W^{1,q} (\Om )$ to $W^{1,q} (\R^n)$  (see \cite[Theorem 5.28]{AdamsFournier}).

\begin{lemma}[Sobolev inequality in Lorentz spaces]\label{Sobolev}
Let $\Om$ be a bounded Lipschitz domain in $\R^2$ and $1< q<2$, there
are constants $C_1=C_1(q,\Om)>0$ and $C_2=C_2(q)$ such that
\EQN{
\|u \|_{L^{q^* , q}(\Om )} &\le C_1  \| u \|_{W^{1,q}(\Om)},\quad \forall u \in W^{1,q} (\Om ),\\
\|u \|_{L^{q^* , q}(\Om )} &\le C_2  \| u \|_{W^{1,q}_0(\Om)},\quad \forall u \in W^{1,q}_0 (\Om ).
}
\end{lemma}
This lemma will allow us to estimate the drift term later. The constant $C_2$ does not depend on $\Om$ because we can take zero extension outside of $\Om$ when $u \in W^{1,q}_0 (\Om )$. This uniformness will be useful for our boundary estimate.

Next we prove $L^2$ estimates for $q$-weak solutions, $q>2$.
\begin{lemma}[Energy]\label{energyest}
  Let $\Om$ be a bounded Lipschitz domain in $\R^2$, $2<q<\infty$, and $u$ be a $q$-weak solution of \eqref{drifteq} for $\GG\in L^q(\Om)$. Then we have the following energy estimate
  \EQ{\label{eq2.2}
  \|u\|_{W^{1,2}(\Om)}\le c\|\GG\|_{L^q(\Om)},
  }
with $c$ depending on $|\Om|$ only.
\end{lemma}
\begin{proof}
  Since $\div b=0$ and $u\in W^{1,q}(\Om)$ we can use $u$ as test function in the weak form \eqref{weak-soln} of \eqref{drifteq}, and justify
 $\int_\Om b \otimes u : \nb u =-\frac12\int_\Om b\cdot \nb |u|^2=  0$. From $\int|\nb u|^2 = - \int \GG: \nb u $ we get $\norm{\nb u}_{L^2(\Om)} \le \norm{\GG}_{L^2(\Om)}$ and \eqref{eq2.2} by Sobolev and H\"older inequalities.
\end{proof}

Lastly we formulate a lemma concerning pressure.
\begin{lemma}[Pressure]\label{pressure}
Let $\Om\subset \R^2$ be a bounded Lipschitz domain, and $1<q<\infty$, then for any $f \in W^{-1,q}(\Om)^2$ such that
$
\bka{f,\zeta}=0$ for all $\zeta \in W^{1,q'}_{0,\si}(\Om)$,
there exists a unique $p \in L^q_0(\Om)$, such that $\norm{p}_{q} \le C \| f\|_{-1,q}$ and
\[
\bka{f,\zeta}=\int_\Om p \div \zeta, \quad \forall \zeta \in W^{1,q'}_{0}(\Om)^2.
\]
\end{lemma}
This is a special case of \cite[Theorem III.5.3]{Galdi} for bounded domains.

\section{Existence of weak solutions}\label{sec3}
In this section we will prove existence of weak solutions to the perturbed Stokes system \eqref{drifteq} for $b\in L^{2,\infty}$, $\div b=0$,
\EQ{ \label{eq3.1}
-\Delta u+ b\cdot\nabla u + \nb \pi =f,\quad \div u=0, \quad u|_{\pd \Om}=0,
}
in $\Om \subset \R^2$.
Recall that a weak solution $u$ belongs to $ W^{1,2}_{0,\sigma}(\Omega)$ and satisfies a weak form of \eqref{drifteq} similar to \eqref{weak-soln}.
\begin{theorem}[Weak solutions  for weak $L^2$ drift]\label{energy}
Let $\Om$ be a bounded Lipschitz domain in $\R^2$.
Assume $b\in L^{2,\infty}(\Omega)^2$ and $\div b=0$.
Then for any $f \in W^{-1,2}(\Omega)^2$,
there exists a weak solution $u\in W^{1,2}_{0,\sigma}(\Omega)$
of \eqref{drifteq}. Moreover, there is $\pi \in L^2(\Om)$ so that $(u,\pi)$ solves \eqref{drifteq} in distributional sense, $\int_\Om \pi=0$, and
\EQ{\label{eq6.1}
\|u\|_{ W^{1,2}(\Omega)}+\|\pi\|_{L^2(\Omega)}\le C\|f\|_{W^{-1,2}(\Om)},
}
for some constant $C=C(\Om)$ independent of $b$.
\end{theorem}

It is important to note that Theorem \ref{energy} does not claim uniqueness, and estimate \eqref{eq6.1} may fail for another solution. %
In fact, in scalar PDE case \eqref{scalar.eq}, there are results that shows non-uniqueness for weak solutions if drift is weaker than $L^2$, see for example \cite[\S2]{Zhikov} for $b\in L^{3/2-}$. Later in Section 6 we will prove uniqueness of $q$-weak solution  for $b\in L^{2,\infty}$ by a duality argument for the case $q<2$, from which follows the uniqueness in case $q=2$.

We will use the following approximation theorem concerning Lipschitz domains in $\R^n$, $n\ge2$. See
Verchota \cite[Theorem 1.12]{Verchota_2} or Kwon \cite[Theorem 2.1]{Kwon} for the statement, and Verchota \cite[Appendix]{Verchota_1} for the proof.
A \emph{coordinate cylinder} $Z$ for $\pd \Om$ is a right circular cylinder with center $Q\in \pd \Om$ such that in a suitable coordinate system (after rotation and translation) with center $Q=(0,0)$, we have $Z=\{(x',x_n)\in\R^n: |x'|<r, |x_n|<H\}$ for some $r,H>0$ and
$ \Om \cap Z=\{(x',x_n)\in\R^n: |x'|<r,  \gamma(x')<x_n <H\}$ for some function $\gamma$ on $\R^{n-1}$ with compact support and $|\ga|< H$. We denote the scaled cylinder $kZ = \{ Q+k(x-Q): x \in Z\}$ for $k>0$.

\begin{theorem}[Approximations of a Lipschitz domain]\label{Lip_domains}
  Let $\Om$ be a bounded Lipschitz domain in $\Bbb R^n, ~n\ge 2$. There exist a finite covering $\mathcal Z = \{ (Z_i,\gamma_i)\}$ of $\pd \Om$ of coordinate cylinders, and a
  sequence of $C^{\infty}$-domains $\{\Om_j\}_j $ with $\bar{\Om}_j\subset\Om $ such that
  with $M=\max_{i} \norm{\nb \gamma_i}_{L^\infty(\R^{n-1})}$,
   for any $(Z,\gamma)\in \mathcal Z$, $Z^*= 100(M+1)Z$
   is still a coordinate cylinder,
  $ \Om \cap Z^*$ is given as $x_n > \gamma(x')$ in a suitable coordinate system, and
   for each $j$,
 $\partial \Om_j\cap Z^*$ is given as a graph of a $C^{\infty}$ function $\gamma_j$, such that $\gamma_j\rightarrow\gamma$ uniformly, $\|\nabla\gamma_j\|_{L^{\infty}(\R^{n-1})}\le \|\nabla\gamma\|_{L^{\infty}(\R^{n-1})}$, and $\nabla\gamma_j\rightarrow\nabla\gamma$ pointwise a.e..
\end{theorem}

We do not need the full power of the theorem. We only need a sequence of subdomains $\Om_n\nearrow \Om$ with uniformly bounded Lipschitz constants and a uniform radius $R_0$ (see \eqref{R0.def}) for balls centered at boundary so that the boundary is a Lipschitz graph in any of these balls.

\medskip

Now we prove the existence theorem for drift in $L^{2,\infty}$.
\begin{proof}[Proof of Theorem \ref{energy}]
  By Theorem \ref{Lip_domains} we have a sequence of smooth subdomains $\Om_n\Subset\Om$ and $\cup_n\Om_n=\Om$. Next we approximate the drift using standard mollification with kernel $\omega\in C^{\infty}_c(B_1)$, $\omega\ge 0,~\int_{\Bbb R^2}\omega=1$. Let $\om_\e(x)=\e^{-2} \om(\e^{-1}x)$. We choose $\varepsilon_n\rightarrow 0^+,\varepsilon_n\le \dist\{\Omega_n,\partial\Om\}$ and extend $b(x)=0$ for $x\not\in \Om$, then we denote
  \EQ{
  b_n(x)=\int_{\Om}\omega_{\varepsilon_n}(x-y)b(y)\,dy,~\forall x\in \Omega.
  }
  From the properties of mollification we have the following
  \EQ{\label{b_n_weak}
  \|b_n-b\|_{L^r(\Omega)}\rightarrow 0,~\forall 1<r<2,~\text{and}~ \|b_n\|_{L^{2,\infty}(\Om)}\le c\|b\|_{L^{2,\infty}(\Om)}.
  }
  Moreover on the sub-domain $\Om_n$ we have $\div b_n=0$ and since $b_n$ is a smooth function in $\Om_n$ there exists a unique weak solution $u_n,\pi_n$ in $\Om_n$ that solves
  \EQ{\label{ap_seq}
  -\Delta u_n+(b_n\cdot\nabla)u_n+\nabla\pi_n=f,\quad\div u_n=0 \quad\text{in}~\Om_n,\quad u_n|_{\partial\Om_n}=0.
  }
 We can extend $u_n$ by zero to $\Om \setminus \Om_n$.
  Since $b_n$ is smooth we can multiply \eqref{ap_seq} by $u_n$ and integrate to get that these solutions also satisfy $W^{1,2}$ bound
  \EQ{\label{w_12_weakbound}
  \|u_n\|_{W^{1,2}(\Om)}\le C\|f\|_{W^{-1,2}(\Om)}.
  }
  The final step is to use this a priori bound to construct a weak $W^{1,2}$ limit $u_n\rightharpoonup u$ and we need to check that the limit will satisfy the equation in the sense of distributions. Let $\xi\in C^{\infty}_c(\Om)$ such that $\div \xi=0$. From the properties of our domains $\Om_n$ we can find $N$ such that for any $n\ge N$, $\supp \xi\subset\Om_n$, which makes the function $\xi$ an admissible test function in \eqref{ap_seq}. Combining it with \eqref{b_n_weak} we can take limits of the weak form of \eqref{ap_seq} and get that $u$ is a weak solution of \eqref{drifteq}.
\end{proof}

\section{Wolf's local pressure decomposition}\label{sec4}
For internal and boundary estimates we will need to deal with pressure terms, for that we will use Wolf's local pressure projection \cite{Wolf15, Wolf17,JWZ}.
Consider the unperturbed Stokes system in a domain $V\subset \R^n$, $n\ge2$,
\EQ{\label{Stokes}
-\Delta v +\nabla \pi=f,\quad \div v=0, \quad v|_{\pd V}=0.
}
We use the same definition of weak solution pairs \eqref{soln-pair} and $p$-weak solutions \eqref{weak-soln} of \eqref{Stokes} in the same way as for \eqref{drifteq}. If there is a unique weak solution pair $(v,\pi)$ with $v \in W^{1,p}_{0,\si}(V)$ and $\pi \in L^p_0(V)$ for each given $f \in W^{-1,p}(\Om)^n$, then
the map
\[
\mathcal{W}_{p,V}(f)=\nabla\pi
\]
is the local pressure projection in $V$.

\begin{theorem}\label{theoremL2wolf}
Assume $p>1$ and $V\subset \R^n$, $n\ge2$, is a bounded Lipschitz domain with Lipschitz constant less than a sufficiently small  $L_0(p,n)>0$.
Then the problem \eqref{Stokes} has a unique weak solution pair $(v,\pi)$ with $v \in W^{1,p}_{0,\si}(V)$ and $\pi \in L^p_0(V)$ for given $f=\div \FF$, $ \FF \in L^p(\Om)^{n\times n}$, and the local pressure projection in $V$,
\[
\mathcal{W}_{p,V}(\div \FF)=\nabla\pi,
\]
satisfies the following inequality
\EQ{\label{pWolfbound}
\|\pi\|_{L^p(V)}\le c_1\|\FF\|_{L^p(V)},
}
where $c_1=c_1(n,V,p)$.
\end{theorem}

This result is well known for a more regular domain, and originated from the classical work of Cattabriga \cite{Cattabriga} for $C^2$-domains in $\R^3$, but for Lipschitz domains we refer to Galdi, Simader, and Sohr \cite{GSS94}, which also considers nonzero $\div v$ and boundary value of $v$.

Compared to \cite[Theorem 4.2]{CT1}, the exponent $p$ is arbitrary. However,
the constant $c_1$ in \eqref{pWolfbound} depends on the domain $V$ and its dependence on $V$ is implicit. We also need to assume small Lipschitz constant.
In contrast, in Theorem 4.2 in \cite{CT1}, the exponent is limited to $p=2$, but it has $L^2$ estimates with uniform constant over any Lipschitz domain with Lipschitz constant less than $L_0<1/2$.

In the current paper for dimension 2, we cannot use the local pressure projection $\mathcal{W}_{p,V}$ for $p=2$, and we need to choose $1<p<2$, as explained in Section \ref{sec1} due to $b\otimes u\in L^p$ for $p<2$ and not $p=2$. Therefore, we will fix $p=4/3$. To retain a uniform constant, we will only apply $\mathcal{W}_{p,V}$ to two choices of $V$: the unit disk and a smooth domain in between $B_{5/8}^+$ and $B_{1}^+$, given by
the following geometrical lemma.

\begin{lemma}\label{geometry}
For any $0<\rho<R$, there is a smooth domain $V\subset \R^2$ such that
\[
B_\rho ^+ \subset V \subset B_R^+.
\]

\end{lemma}
This geometric property is well known, so we will leave its proof in the Appendix.

\begin{remark}\label{Wolfpdecom}
We will choose $p=4/3$ and apply Wolf's local pressure projection in $V$ to equation  \eqref{drifteq} to get
\EQ{
\nb \pi = \cW_{4/3,V}(\nb \pi ) =\cW_{4/3,V}[\Delta u-\div(b\otimes u)+\div \GG].
}
Thus we have a local pressure decomposition $\pi-(\pi)_{V} =\pi_1+\pi_2+\pi_3$ in $V$, where $\pi_1,\pi_2,\pi_3\in L^{4/3}_0(V)$ are given by
\EQ{
\nabla\pi_1=\mathcal{W}_{4/3,V}(\Delta u),\quad \nabla\pi_2=\mathcal{W}_{4/3,V}(-\div(b\otimes u)),\quad \nabla\pi_3=\mathcal{W}_{4/3,V}(\div \GG).
}
In the internal case we will only have three terms in decomposition, but in the boundary case there will be also error terms due to boundary stretching when $V$ is between  $B_{5/8}^+$ and $B_{1}^+$.  By Theorem  \ref{theoremL2wolf}, we have for $c=c_1(2,V,4/3)$
\EQ{\label{eq4.14}
\norm{\pi-(\pi)_{V}}_{L^{4/3}(V)} \le c\bke{ \norm{\nb u}_{L^{4/3}(V)}
+ \norm{b\otimes u}_{L^{4/3}(V)} + \norm{\GG}_{L^{4/3}(V)}}.
}
\end{remark}

\section{A priori bound for large drift and $p$ close to $2$}\label{sec5}

In this section we prove the key a priori bound for higher integrability of weak solutions
assuming higher integrability of the drift.

\begin{theorem}\label{Capriori}
Let $\Om$ be a bounded Lipschitz domain in $\R^2$ with Lipschitz constant $L \le L_0$, where $0<L_0<1$ is a sufficiently small constant. Assume $b\in L^{2,\infty}(\Omega)^2$ and $\div b=0$.
Assume further $b \in L^{2+\delta}$ for some $\delta>0$.
Then there exists $p_0(\|b\|_{L^{2,\infty}})>2$ such that for any $2<p<p_0$,
any weak solution $u\in W^{1,2}_{0,\si}(\Omega)$  of \eqref{drifteq} with force $\GG\in L^p(\Om)^{2\times 2}$
is in $ W^{1,p}(\Omega)$ and
\EQ{\label{eq5.1}
\|u\|_{W^{1,p}(\Omega)}\le c\|u\|_{W^{1,2}(\Omega)}+ c\|\GG\|_{L^p(\Omega)}
}
for some constant $c(\Omega,\|b\|_{L^{2,\infty}})>0$.
\end{theorem}

\emph{Comments on Theorem \ref{Capriori}:}
\EN{
\item Estimate \eqref{eq5.1} can be reduced to
\EQ{\label{eq5.2}
\|u\|_{W^{1,p}(\Omega)}\le  c\|\GG\|_{L^p(\Omega)}
}
if $u$ is a weak solution constructed in Theorem \ref{energy}. For a general weak solution, we need to keep the term $c\|u\|_{W^{1,2}(\Omega)}$ on the right side of \eqref{eq5.1} for now. It can be removed after we have proved Theorem \ref{main1}.

\item Although we assume $b \in L^{2+\delta}$, the constants only depend on $\|b\|_{L^{2,\infty}(\Omega)}$.
The assumption $b \in L^{2+\delta}$ can be replaced by $\nb u \in L^{2+\delta}$. Any of these assumptions makes the first integral in
\EQ{\label{5.1}
\int_\Om b \otimes u : \nb( u \eta^4) = -\int_\Om b | u|^2\cdot 2\eta^3\nb \eta
}
 integrable. Once integrable, the equality can be shown,
and the second integral is defined without the extra assumption.

\item The dependence of  the constants $p_0$ and $c$ on $b$ and $\Om$ is only on $\|b\|_{L^{2,\infty}}$, diam\,$\Om$, and the radius $R_0$ in \eqref{R0.def}
so that $\pd\Om \cap B_{2r}(x_0)$ is a Lipschitz graph for any $x_0\in\pd\Om$ and $r<R_0$.

\item The smallness of $L_0\le1$ is for boundary stretching for local pressure estimate, and is explicitly specified after \eqref{change-var} and after \eqref{pressure_bound_est}. It is not needed for the scalar equation case in Theorem \ref{th1.2}.
}

\begin{proof}
Similar to \cite{CT1} we will use the Gehring's lemma approach \cite{Gehring, Gia}. Denote $\Om_{x,R}=\Om \cap B(x,R)$ for $x\in \overline \Om$ and $R>0$. Since $\pd \Om$ is Lipschitz and compact, there
is a radius $R_0\in (0,1]$ such that for every $x_0 \in\pd \Om$,
\EQ{ \label{R0.def}
\Om_{x_0,2R_0}= \{ (x',x_n) \in B_{2R_0}(0) \subset \R^2,\ x_n>\ga(x')\},
}
after suitable coordinate rotation and translation, where $\ga(x')$ is a Lipschitz function defined for $x'\in B_{2R_0}'(0) \subset \R$ with Lipschitz constant $L$ and $\ga(0)=0$.

To prove the a priori estimate \eqref{eq5.1}, we first consider the interior case with $x_0\in\Omega$ and $B_{\rho}=B_{\rho}(x_0)\subset\Omega$. We rescale it to the unit ball $B$ and denote $u_{\rho}(x)=\rho u(\rho x+x_0),b_{\rho}(x)=\rho b(\rho x+x_0),\pi_{\rho}(x)=\rho^2\pi(\rho x+x_0)$ and $\GG_{\rho}(x)=\rho^2 \GG(\rho x+x_0)$. Let $\bar{u}_{\rho}=u_{\rho}-k$ where $k$ is a real constant to be choose later.  By Remark \ref{Wolfpdecom}, we can decompose $\pi_{\rho} - (\pi_{\rho})_{B}=\pi_1+\pi_2+\pi_3$ in $B$, where $\pi_1,\pi_2,\pi_3$ are given by
\EQ{
\nabla\pi_1=\mathcal{W}_{4/3,B}(\Delta \bar{u}_{\rho}),\quad \nabla\pi_2=\mathcal{W}_{4/3,B}(-\div(b_{\rho}\otimes \bar{u}_{\rho})),\quad \nabla\pi_3=\mathcal{W}_{4/3,B}(\div \GG_{\rho}),
}
with $\int_B \pi_i=0$.
Therefore, by Theorem \ref{theoremL2wolf} we have the following pressure bound
\EQ{\label{wolfbound}
\|\pi_1\|_{4/3,B}\le c_1\|\nabla \bar{u}_{\rho}\|_{4/3,B},\quad\|\pi_2\|_{4/3,B}\le c_1 \|b_{\rho}\bar{u}_{\rho}\|_{4/3,B},\quad \|\pi_3\|_{4/3,B}\le c_1\|\GG_{\rho}\|_{4/3,B},
}
where $c_1$ is a global constant.
For simplicity of notations let us drop the subscript $\rho$ for the further calculations, until we need to scale back to ball $B_{\rho}(x_0)$.

Take a smooth cutoff function $\eta$ on $B$ with $\eta=1$ in $B_{\frac{1}{2}}$. Use the test function $\zeta=\eta^4 \bar{u}$ in the weak form \eqref{soln-pair} with $(u,\pi)$ replaced by $(\bar u, \pi - (\pi)_B)$ to get the energy estimate (using \eqref{5.1} due to $b \in L^{2+\de}$)
\EQS{\label{3.1}
\int_{B} |\nabla (\eta^2\bar{u})| ^2
\le &\int_{B} |\bar{u}|^2(|\nabla \eta^2|^2) +\int_{B} |\bar{u}|^2|b||\nabla\eta|\eta^3
\\
&+\int_{B} (|\pi_1|+|\pi_2|+|\pi_3|)|\bar{u}||\nabla\eta|\eta^3+\bka{\div \GG_{\rho}, \bar{u}\eta^4}_B.
}

We will estimate each of the terms on the right hand side separately, using Lemmas \ref{Holder} and \ref{Sobolev} in ball $B$,
\EQS{
\int_{B} |\bar{u}|^2|b||\nabla \eta|\eta^3
&\le \|b\|_{L^{2,\infty}(B)}\|\bar{u}^2\eta^3\|_{L^{2,1}(B)}\le \|b\|_{L^{2,\infty}(B)}\|\bar{u}\eta^{3/2}\|_{L^{4,2}(B)}^2
\\
&\le C \|b\|_{L^{2,\infty}(B)}\|\bar{u}\eta^{3/2}\|^2_{W^{1,4/3}(B)}.
}
For pressure terms we use \eqref{wolfbound} and get the following
\EQN{
\int_{B}|\pi_1||\bar{u}||\nabla\eta|\eta^3 &\le \|\pi_1\|_{4/3,B}\|\eta^2 \bar{u}\|_{L^4(B)}
\\
&\le
c\|\nabla \bar{u}\|_{4/3,B}\|\eta^2 \bar{u}\|_{W^{1,4/3}(B)}\le c\|\nabla \bar{u}\|_{4/3,B}^2+c\|\eta \bar{u}\|^2_{4/3,B},
}
\EQN{
\int_{B}|\pi_2||\bar{u}||\nabla\eta|\eta^3 &\le \|\pi_2\|_{L^{4/3}(B)}\|\eta^2\bar{u}\|_{L^4(B)}\le c\|\bar{u}\|_{L^{4,4/3}(B)}\|b\|_{L^{2,\infty}(B)}\|\eta^2\bar{u}\|_{L^{4}(B)}
\\
&\le c(b)\|\nabla \bar{u}\|_{4/3,B}^2+c(b)\| \bar{u}\|^2_{4/3,B}
}
\EQN{
\int_{B}|\pi_3||\bar{u}||\nabla\eta|\eta^3 &\le \|\pi_3\|_{4/3,B}\|\eta^2 \bar{u}\|_{L^4(B)}\le  c\|\GG\|_{2,B}\|\eta^2 \bar{u}\|_{W^{1,4/3}(B)}
\\
&\le \frac 14 \|\GG\|_{2,B}^2+c\|\eta^2 \bar{u}\|_{W^{1,4/3}(B)}^2.
}
Lastly we estimate the term containing $\GG$,
\EQN{
\bka{\div \GG_{\rho}, \bar{u}\eta^4}_B
&=-\int_{B} \GG :\nabla(\bar{u}\eta^4)=-\int_{B} \GG: \nabla(\eta^2\bar{u})\eta^2-2\int_{B} \GG :\bar{u}\nabla\eta\eta^3
\\
&\le \frac{1}{16}\|\nabla(\eta^2\bar{u})\|^2_{2,B}+c\|\eta^2 \bar{u}\|_{2,B}^2+ \frac 14\|\GG\|_{2,B}^2.
}
Combining these estimates with \eqref{3.1} we get the following
\EQ{
\int |\nabla (\eta^2\bar{u})|^2\le c(b)\|\bar{u}\|^2_{W^{1,4/3}(B)}+c(b)(\|\bar{u}\|_{2,B}^2+\|\bar{u}\|_{4/3,B}^2)+\|\GG\|_{2,B}^2,
}
Here $c(b)$ depends on $b$ only through $\norm{b}_{L^{2,\infty}(\Om)}$. Finally, since $\eta=1$ on $B_{1/2}$ we get that
\EQ{
\int_{B_{1/2}} |\nabla \bar{u}|^2\le c(b)\|\nabla \bar{u}\|^2_{4/3,B}+c(b)\|\bar{u}\|_{2,B}^2+\|\GG\|_{2,B}^2.
}
We now choose $k=(u_{\rho})_{B}$ as the constant in $\bar{u}=u_{\rho}-k$ and apply Poincar\'e inequality to get
\EQ{
\int_{B_{1/2}} |\nabla u_{\rho}|^2\le c(b)\|\nabla u_{\rho}\|^2_{4/3,B}+\|\GG\|_{2,B}^2.
}
Lastly we scale back to $B_{\rho}(x_0)$ and get that
\EQ{\label{5.12}
\frac{1}{|B_{\rho/2}|}\int_{B_{\rho/2}} |\nabla u|^2
\le c(b) \Big(\frac{1}{|B_{\rho}|}\int_{B_{\rho}}|\nabla u|^{4/3}\Big)^{3/2}+\frac{c}{|B_{\rho}|}\|\GG\|_{2,B_{\rho}}^2.
}

\medskip

Next we consider the boundary case, $\Omega_{\rho}=\Omega\cap B_{\rho}(x_0)$ with $x_0\in\pd\Om$, and $0<\rho<R_0$. The significant difference from the internal case is the pressure estimate: Unlike the internal case we can not use Wolf decomposition estimate in a uniform ball, therefore we need to track the dependence of constant $c_1$ when we apply Theorem \ref{theoremL2wolf} to get estimates similar to \eqref{pressure_bound_est}
 and make sure it is uniform across all  $\Omega_{\rho}(x_0)$.
 Rescale the domain as the interior case so that $x_0=0$ and $\rho=1$.
 We will use change of coordinates
 \EQ{\label{change-var}
 y_1=x_1, \quad y_2=x_2-\gamma(x_1),
 }
 to straighten the boundary and map $\Om_{5/8}\subset\Om_1$ into domains $\Om_{5/8}'\subset\Om_1'$ with flat boundary on $y_2=0$. Since $\gamma$ was a Lipschitz curve with sufficiently small Lipschitz constant we can assume that $\Om_{5/8}'\subset B_{6/8}^+\subset B_{7/8}^+\subset\Om_1'$.
  (This is the first place that specifies the smallness of the Lipschitz constant.) By Lemma \ref{geometry} we can fix a smooth $V$ such that $B_{6/8}^+\subset V\subset B_{7/8}^+$. After rewriting the equation \eqref{drifteq} in new variables we get the following system (with all derivatives in $y$) in $\Om_1'$:
\EQS{\label{changed_var}
-\Delta_y u+\div (\gamma' e_1\otimes\partial_{y_2}u)+
\gamma' \pd_1\pd_2 u - (\gamma')^2\partial_{y_2}^2u+
\\
+\div_y(b\otimes u)- \gamma'\partial_{y_2}(b\otimes u)_1+\nabla_y \pi-\gamma'\partial_{y_2}\pi e_1&=\div\GG-\gamma'\partial_2 \GG_1 ,
\\
\div_y u-\gamma'\partial_2 u_1&=0.
}
Above, $(b\otimes u)_1=(b_1u_1,b_1u_2)$ and $\GG_1=(G_{11},G_{12})$ are their first rows. In the second term, $\ga'$ is hidden behind $\div$ to avoid $\ga''$ which we do not assume any bound.

We apply Wolf's local pressure projection $\mathcal{W}_{4/3,V}$ on \eqref{changed_var} in the smooth domain $V$ and get that  $\pi-(\pi)_{V}=\pi_1+\pi_2+\pi_3+\pi_4$ where
\EQS{\label{pressure_decomp_b}
\nabla\pi_1&=\mathcal{W}_{4/3,V}\bkt{\Delta_y u-\div (\gamma' e_1\otimes\partial_{y_2}u)-\gamma' \pd_1\pd_2 u +(\gamma')^2\partial_{y_2}^2u}
\\
\nabla\pi_2&=\mathcal{W}_{4/3,V}\bkt{-\div_y(b\otimes u)+\gamma'\partial_{y_2}(b\otimes u)_1}
\\
 \nabla\pi_3&=\mathcal{W}_{4/3,V}(\gamma' \partial_{y_2}\pi e_1)
 \\
 \nabla\pi_4&= \mathcal{W}_{4/3,V}(\div \GG-\gamma'\partial_2 \GG_1).
}
From Theorem \ref{theoremL2wolf} the local pressure projection  $\mathcal{W}_{4/3,V}$ is well defined, and we have the bound \eqref{pWolfbound} with global constant $c$. Also using $\ga'\pd_2 = \pd_2 \ga'$, and Lemmas \ref{Holder} and \ref{Sobolev} in $V$,
\EQS{\label{pressure_bound_est}
\|\pi_1\|_{4/3,V}&\le c(1+\|\gamma'\|_{\infty}+\|\gamma'\|_{\infty}^2)\|\nabla u\|_{4/3,V}
\\
\|\pi_2\|_{4/3,V}&\le  c(1+\|\gamma'\|_{\infty})\|u\otimes b\|_{4/3,V}
\\
\|\pi_3\|_{4/3,V}&\le  c\|\gamma'\|_{\infty}\|\pi\|_{4/3,V}\le c\|\gamma'\|_{\infty}\|\pi_3\|_{4/3,V}+c\|\gamma'\|_{\infty}\|\pi_1+\pi_2+\pi_4\|_{4/3,V}
\\
\|\pi_4\|_{4/3,V}&\le  c(1+\|\gamma'\|_{\infty})\|\GG\|_{4/3,V}
}
Here we apply our choice of $\gamma$ so that $\|\gamma'\|_{\infty}<1$ and
$ c\|\gamma'\|_{\infty}<1/2$. (This is the second place that specifies the smallness of the Lipschitz constant.) After applying it to \eqref{pressure_bound_est} we get the following
\EQS{\label{pressure_3}
\|\pi_3\|_{4/3,V}\le 2c\|\gamma'\|_{\infty}\|\pi_1+\pi_2+\pi_4\|_{4/3,V}.
}

Let $\Om^\sharp$ be the pre-image of $V$ under the map \eqref{change-var} in $x$-variables.
The pressure component $\pi_i$ can be considered functions on $\Om^\sharp$ and since \eqref{change-var} has determinant 1 and $\|\nabla_y u\|_{4/3,V}\le 2\|\nabla_x u\|_{4/3,\Omega^\sharp}$, by assumption $\|\gamma'\|_{\infty}\le 1$ we have
\EQS{\label{pressure-est}
\norm{\pi-(\pi)_{V}}_{4/3, \Om^\sharp} &\le {\textstyle\sum}_{i=1}^4 \norm{\pi_i}_{4/3, \Om^\sharp}
={\textstyle\sum}_{i=1}^4 \norm{\pi_i}_{4/3, V}
\\
&\le c\|\nabla _y u\|_{4/3,V} +c \|u\otimes b\|_{4/3,V} +c \|\GG\|_{4/3,V}
\\
&\le c(1+ \| b\|_{L^{2,\infty}(V)}) (\|\nabla_y u\|_{4/3,V} +\| u\|_{4/3,V} )+ c \|\GG\|_{4/3,V}
\\
&\le c(1+ \| b\|_{L^{2,\infty}(\Omega)}) (\|\nabla_x u\|_{4/3,\Omega^\sharp} +\| u\|_{4/3,\Omega^\sharp} )+ c \|\GG\|_{4/3,\Omega^\sharp}.
}
This pressure estimate is the only purpose of our change of variables \eqref{change-var}. Noting that $\Om_{5/8}\subset \Om^\sharp\subset \Om_1$, we then proceed the following.

 Let $\eta$ be a cutoff function supported in $B_{5/8}(x_0)$ such that $\eta=1$ on $B_{1/2}(x_0)$. Due to our boundary condition $u|_{\partial\Omega}=0$ we can use function $\eta^4u$ as an admissible test function in \eqref{drifteq}. Therefore we get the following local energy inequality
\EQS{\label{eq5.19}
\int_{\Om_1} |\nabla (\eta^2u)|^2
&\le \int_{\Om_1} |u|^2|\nabla \eta^2|^2+\int_{\Om_1} |u|^2|b||\nabla \eta|\eta^3\\
&\quad +\int_{\Om_1} |\pi-(\pi)_{V}||u||\nabla\eta|\eta^3-\int_{\Om_1} \GG :\nabla(u\eta^4).
}
Using \eqref{pressure-est} to bound the pressure term and similarly to internal case
using Lemmas \ref{Holder} and \ref{Sobolev} in $\Om_1$, with uniform constant $C_2$ in Lemma \ref{Sobolev} using $\eta^4u\in W^{1,2}_0(\Om_1)$, we get
\EQS{\label{5.16}
\int_{\Om_{1/2}} |\nabla u|^2 &\le C(b)(\|\nabla u\|^2_{4/3,\Om_1}+\|u\|_{2,\Om_1}^2)+c\|\GG\|_{2,\Om_1}^2.
}
Notice that here the constant $C(b)$ does not depend on $\Om_1$.
Finally we extend function $u,\GG$ by zero to  $B_1(x_0)\setminus \Omega$. Since $u=0$ on the set $E:=\Om^c\cap B_1(x_0)$ of non zero measure with the ratio $\frac{|E|}{|B_{1}|}$ bounded from below, we can use Poincar\'e inequality to bound $\|u\|_{2,\Om_1}^2 \le c \|\nabla u\|^2_{4/3,\Om_1}$, with a constant independent of $\Om_1$ with the given lower bound of the ratio.
(This version of Poincar\'e inequality follows from Proposition 6.14 of \cite[page 116]{Lieberman},  with $f=\al 1_E$, %
$\al=\frac {|B_1|}{|E|} $, and $\fint_{B_1} f=1$. It gives $\norm{u}_{4/3,B_1} \le C\al  \norm{\nabla u}_{4/3,B_1} $.
We also have
$\norm{u}_{4,B_1} \le C \norm{u}_{W^{1,4/3}(B_1)}$.)

We also re-scale the inequality for arbitrary $B_r(x_0)$ and get \eqref{5.12}.

\medskip
Lastly we need to consider the case of arbitrary $x_0,r$ such that $\Om^c\cap B_{2r}(x_0)\neq\emptyset$. For any such $x_0$ we can find $y_0\in\partial\Om\cap B_{2r}(x_0)$ such that
\EQ{
B_r(x_0)\subset B_{3r}(y_0)\subset B_{6r}(y_0)\subset B_{8r}(x_0).
}
Using \eqref{5.12} for $B_{3r}(y_0)$ (assuming $6r<R_0$) we get
\EQS{\label{5.17}
\frac{1}{|B_r|}\int_{B_r(x_0)} |\nabla u|^2
&\le \frac{C}{|B_{3r}|}\int_{B_{3r}(y_0)} |\nabla u|^2
\\
&\le c(b) \Big(\frac{1}{|B_{6r}|}\int_{B_{6r}(y_0)}|\nabla u|^{\frac{4}{3}}\Big)^{\frac{3}{2}}+\frac{C}{|B_{6r}|}\|\GG\|_{L^2(B_{6r}(y_0))}^2
\\
& \le c(b) \Big(\frac{1}{|B_{8r}|}\int_{B_{8r}(x_0)}|\nabla u|^{\frac{4}{3}}\Big)^{\frac{3}{2}}+\frac{C}{|B_{8r}|}\|\GG\|_{L^2(B_{8r}(x_0))}^2.
}

For Gehring's lemma we combine the internal estimate \eqref{5.12} and the boundary estimate \eqref{5.17} and conclude that
\EQ{\label{5.19}
\frac{1}{|B_r|}\int_{B_r(x_0)} |\nabla u|^2 \le
c(b) \Big(\frac{1}{|B_{8r}|}\int_{B_{8r}(x_0)}|\nabla u|^{\frac{4}{3}}\Big)^{\frac{3}{2}}+\frac{C}{|B_{8r}|}\|\GG\|_{L^2(B_{8r}(x_0))}^2,
}
for any $x_0$ in a big cube containing $\overline \Om$, and any $r\le R_0/6$.
We can now apply Proposition 1.1 of Giaquinta \cite[page 122]{Gia} (also see \cite[Proposition 3.7]{DongKim}) to get that $\nb u \in L^p_\loc$ for $p \in (2,p_0)$ for some $p_0 >2$, and
\EQ{\label{0619}
\bke{\fint_{B_r(x_0)} |\nabla u|^p}^{\frac 1p} \le
c_1\bke{\fint_{B_{8r}(x_0)} |\nabla u|^2}^{\frac 12} + c_1\bke{\fint_{B_{8r}(x_0)} |\GG|^p}^{\frac 1p},
}
for all $x_0 \in \overline \Om$ and $r\le R_0/6$, with
$p_0$ and $c_1$ depending only on $c(b)$. Summing \eqref{0619} over a finite cover of $\overline \Om$ of balls of radius $R_0/6$,
we get \eqref{eq5.1}.
\end{proof}

\begin{remark}\label{remark5.2}
 Tracking the proof, the constant $p_0$ only depends on  $\norm{b}_{L^{2,\infty}}$, while the constant $c$ in \eqref{eq5.1} only depends on $\norm{b}_{L^{2,\infty}}$, $R_0$, and diam $\Om$, and has no dependence on $\delta$. This is typical for a priori estimates since we assume slightly higher initial regularity of drift which is not used in the final estimate itself. In the next section we will be applying this lemma to sub-domain sequence $\Om_n\subset\Om$, this allows us to have constant being uniform in $n$.
\end{remark}

\section{Proof of Theorem \ref{main1}}\label{sec6}
In this last section we prove the main Theorem \ref{main1} on the unique existence of
weak solution pair $u\in W^{1,q}_{0,\si}(\Omega)$ and $\pi\in L^q_0(\Omega)$ of \eqref{drifteq}
for $q$ sufficiently close to $2$, when $b\in L^{2,\infty}(\Om)$.

\begin{proof}
The existence in the case $q=2$ is already given by Theorem \ref{energy}.
We start with the case $q>2$. The solutions will be constructed with approximation method and a-priori estimates, and we will approximate the Lipschitz domain $\Om$ by sub-domains. By Theorem \ref{Lip_domains}, there exists a sequence $\Om_n$ of Lipschitz sub-domains such that
\EQ{
\Om_n\subset \Om_{n+1} \quad (n \in \NN), \quad\textstyle \bigcup_{n\in\NN} \Om_n=\Om.
}
We also get that Lipschitz constants of $\Om_n$ are uniformly bound by the Lipschitz constant $L$ of $\Om$, and the radius bound $R_0$ specified in \eqref{R0.def} in the proof of Theorem \ref{Capriori} can be taken the same for all $n$.
Choose $\varepsilon_n\rightarrow 0^+$ such that $\varepsilon_n<\dist\{\Om_n,\partial\Om\}$.

Next we mollify the drift $b$: Extend $b(x)=0$ for $x\not \in\Om$ and for any $n$ we take
  \EQ{
  b_n(x)=b_{\varepsilon_n}(x)=\int_{\Om}\omega_{\varepsilon_n}(x-y)b(y)~dy,~x\in\Om.
  }
  Here $\omega_{\varepsilon_n}=\varepsilon_n^{-2}\omega(x/\varepsilon_n)$, and $\omega$ is a standard mollification kernel
  \EQ{
  \omega\in C_c^{\infty}(B_1),\quad\omega\ge0,\quad\int_{\Bbb R^2}\omega=1.
  }
  From the properties of mollification we know that $\div b_n=0$ on $\Om_n$. We also have the following convergence
  \EQN{
  \|b_n\|_{L^{2,\infty}(\Om)}\le c \|b\|_{L^{2,\infty}(\Om)},\quad
  \|b_n-b\|_{L^r(\Om)}\underset{n\rightarrow 0}{\longrightarrow}\infty,~\forall 1<r<2.
  }
   Lastly we also approximate right hand side $\GG$ by $\GG_n\in C_0^{\infty}(\Om_n),~\|\GG_n-\GG\|_{L^q(\Om)}\rightarrow 0$. Since $b_n,\GG_n$ are smooth functions  by Theorem \ref{energy}
  there exists weak solution $(u_n,\pi_n)\in W^{1.2}_0(\Om_n)\times L^2(\Om_n)$ of the following system
  \EQ{\label{aprox_seq}
  -\Delta u_n+(b_n\cdot\nabla)u_n+\nabla\pi_n=\div \GG_n,~\div u_n=0 ~\text{in}~\Om_n,\quad u_n|_{\partial\Om_n}=0,
  }
with $\norm{\nb u_n}_{L^{2}(\Om_n)} \le \norm{\GG_n}_{L^2(\Om_n)}$  and constant 1.
  By Theorem  \ref{Capriori} we know that $u_n\in W^{1,q}_0(\Om_n)$ for $2<q<p_0$, $p_0=p_0(c\| b\|_{L^{2,\infty}(\Om)})>0$ uniform in $n$, and we can extend $u_n$ to the whole $\Om$ by zero. From convergence of $\GG_n\rightarrow \GG$, Lemma \ref{energyest}, Theorem \ref{Capriori} and Remark \ref{remark5.2} we also get that $u_n$ are uniformly bounded in $W^{1,q}(\Om)$,
\EQ{\label{eq6.5}
\|u_n\|_{W^{1,q}(\Omega)}\le c\|\GG\|_{L^q(\Omega)}.
}
Here it's important to notice that constants $p_0$ and $c$ do not depend on $n$ by Remark \ref{remark5.2} using that the diameter, Lipschitz constant and $R_0$ of $\Om_n$ are uniformly bounded by properties of $\Om$, and norm of $b_n$ is uniformly bounded by norm of $b$.

   Finally, by the uniform bound \eqref{eq6.5}, there exists $u\in W^{1,q}(\Om)$ as a weak limit of $u_n$ and we only need to check that we can pass to the limit in \eqref{aprox_seq}. Let $\xi\in C_{c,\sigma}^{\infty}(\Om)$. We can find $N$ such that $\supp\xi\subset\Om_n$ for all $n>N$, therefore $\xi$ is a suitable test function in \eqref{aprox_seq} and we get the following
  \EQ{\label{6.6}
  \int_{\Om} \nabla u_n:\nabla\xi~dx-\int_{\Om}b_n\otimes u_n: \nabla\xi~dx=-\int_{\Om}\GG_n:\nabla\xi~dx.
  }
  Here we can pass to the limit and get that $u$ is a solution to \eqref{drifteq} in terms of distributions for divergence free test functions. Lastly we get the existence of pressure from Lemma \ref{pressure}.
 The a priori estimate \eqref{main1-eq} follows from \eqref{eq6.5}.

So far for $2<q<p_0$ we have shown the existence of $q$-weak solutions satisfying the a priori bound \eqref{main1-eq}.
For uniqueness without assuming \eqref{main1-eq}, suppose there are two  $q$-weak solutions $u$ and $\tilde u$ of \eqref{drifteq} with the same force $\GG$, $q>2$, which may not satisfy \eqref{main1-eq}. Then the difference $w=u-\tilde u$ is a $q$-weak solution of the homogeneous \eqref{drifteq} with $\GG=0$. Since $\nb w \in L^q$, the integral
\[
\int_\Om b \otimes w : \nb w
\]
is absolutely integrable, and can be shown to be zero. Then the usual energy estimate (i.e., \eqref{weak-soln} with $\zeta \to w$ in $W^{1,q}(\Om)$) gives $\int_\Om |\nb w|^2=0$. Hence $w=0$ and $u = \tilde u$.

\medskip

We now consider the case $q<2$. We will first use a duality argument to
prove a priori estimate in $W^{1,q}$ for $q\in(p_0',2)$. Suppose $u$ is a $q$-weak solution of \eqref{drifteq} with right hand side $\div\GG$, i.e., $u$ satisfies $u \in W^{1,q}_{0\,\si}(\Om)$ and \eqref{weak-soln}
\EQ{\label{6.7}
\int_{\Om} (\nb u - b \otimes u): \nb \zeta   = -\int_\Om \GG : \nb \zeta,\quad
\forall \zeta\in C_{c,\sigma}^{\infty}(\Omega).
}
For any $\FF  \in L^{q'}(\Omega)^{n \times n}$, $2<q'<p_0$, let $v$ be the unique $q'$-weak solution of the dual system \eqref{drifteq_dual} with right hand side $\div\FF$, i.e.,
$v \in W^{1,q'}_{0\,\si}(\Om)$ and
\EQ{\label{6.7b}
\int_{\Om} (\nb v + b \otimes v): \nb \eta   = -\int_\Om \FF : \nb \eta,\quad
\forall \eta\in C_{c,\sigma}^{\infty}(\Omega).
}
It is guaranteed by the first part of this theorem that $v$ uniquely exists and
\EQ{\label{6.8}
\|v\|_{W^{1,q'}(\Om)}\le C(b)\|\FF\|_{L^{q'}(\Om)}.
}
One can justify integration by parts and get from \eqref{6.7b}
\EQ{\label{6.9}
\int_{\Om} (\nb \eta - b \otimes \eta ): \nb  v = -\int_\Om \FF : \nb \eta,\quad
\forall \eta\in C_{c,\sigma}^{\infty}(\Omega).
}

Since $b\otimes u \in L^q(\Om)$ by Sobolev embedding Lemma \ref{Sobolev} in Lorentz spaces, and
$\nabla v\in L^{q'}(\Om)$, the term $b\otimes u:\nabla v $ is integrable. So we can  first use $\eta=u$ as a test function in \eqref{6.9}, and then $\zeta=v$ as a test function in \eqref{6.7} to get that
\EQS{\label{6.11}
\int_{\Om}\FF :\nabla u\,dx&=-\int_{\Om}(\nabla u-b\otimes u):\nabla v \, dx
\\
&=\int_{\Om} \GG:\nabla v\,dx\le \|\GG\|_{L^{q}(\Om)}\|\nabla v\|_{L^{q'}(\Om)}\le C(b)  \|\GG\|_{L^{q}(\Om)}\|\FF\|_{L^{q'}(\Om)}.
}
Here we used inequality \eqref{6.8} for the last step.
Note that we have avoided putting $\eta=u$ in \eqref{6.7b} since $b \otimes v: \nb u$ may not be integrable.
Since the matrix function $\FF\in L^{q'}(\Om)$ in \eqref{6.11} is arbitrary we get that
\EQ{\label{smallparpriori1}
\|u\|_{W^{1,q}(\Om)}\le C(b)  \|\GG\|_{L^{q}(\Om)}.
}

This a priori bound for arbitrary $q$-weak solution for $q\in(p_0',2)$
implies in particular the uniqueness of $W^{1,q}$ weak solutions for $q\in(p_0',2)$. Furthermore, it implies uniqueness in case $q=2$, meanwhile existence and the bound was proven in Theorem \ref{energy}.

Now, after we have proved the a priori estimate, we will prove the existence in a similar way to the case $q>2$. Similarly we approximate our Lipschitz domains, drift $b_n\rightarrow b$ with $b_n=b_{\varepsilon_n}$ and we will also approximate right hand side $\GG_n\in C_0^{\infty}(\Om_n),~\|\GG_n-\GG\|_{L^q(\Om)}\rightarrow 0$. Since $\GG_n, ~b_n$ are smooth there exists a unique weak solution $u_n$ of \eqref{aprox_seq} which we extend to the whole $\Om$. Since the right hand side $\GG_n$ is smooth it will also be in $L^q(\Om)$ and therefore we can apply \eqref{smallparpriori1} to get that
\EQ{
\|u_n\|_{W^{1,q}(\Om)}\le C(b_n)  \|\GG_n\|_{L^{q}(\Om)}\le C(b)  \|\GG\|_{L^{q}(\Om)}.
}
Here we used that $b_n$ and $\GG_n$ are uniformly bounded by $b,\GG$ in the corresponding spaces. From the uniform bound there exists a weak limit $u$ of the sequence $u_n$ in $W^{1,q}(\Om)$. Finally
for any $\xi\in C_{c,\sigma}^{\infty}(\Om)$
we can %
choose $N$ so that $\supp \xi\subset\Om_n$ for $n>N$ and apply a similar argument to the case $q>2$,
where we use $\xi$ as a test function in \eqref{6.6} and pass to the limit $n\rightarrow\infty$. This finishes the proof of the main result.
\end{proof}

\begin{remark} In the proof, we use the domain approximation $\{\Om_n\}_n$ so that we can mollify the drift $b$.
The force approximation $\{\GG_n\}_n$ of $\GG$ is not needed for the case $q>2$, but is needed for $q<2$
to get the existence of the approximation sequence $u_n$. To keep the proof uniform we approximate  $\GG$ in both cases.
\end{remark}

\section{Appendix}\label{sec7}
In this section we provide the proof of Lemma \ref{geometry}.
\begin{proof}
We may assume $\rho=1$ by rescaling.
Let $R_1=\min(2,\frac12(R+1))$.

We start with a smooth function $y= f(x)$ with
\[
f(x)=
e^{- \frac 1{x-1}}, \quad ( 1<x<\infty); \qquad f(x)=0, \quad ( 0 \le x \le 1).
\]
It is strictly increasing for $1<x<\infty$ with range $0<y<1$. Thus for $0<y<1$ we have an inverse function
\[
x=1 + \bke{-\log y}^{-1}, \quad (0<y<1).
\]
In polar coordinates, direct calculation using $\th=\arctan \frac yx$ gives for $x>1$
\[
\frac {d\th}{dx} = \frac{1}{1+(y/x)^2} \frac d{dx} \bke{ \frac 1x \,e^{-\frac1{x-1}}}
= \frac{-x^2+3x-1}{(x^2+y^2)(x-1)^2} \, e^{-\frac1{x-1}},
\]
which is positive for $1<x<\frac12(3+\sqrt 5)\approx 2.6180$. At $x=2.5$ we have $y=e^{-2/3}\approx0.5134$, $\th\approx 0.20255$, and $r\approx
2.5522$. Thus for $0<\th\le 0.2$, the curve $y=f(x)$ is described by
\[
r=g(\th), \quad (0<\th\le 0.2).
\]

Let $g_1(\th) = \min (g(\th),R_1)$. Let $\th_0>0$ be the smallest angle such that $g(\th)=R_1$. We have $0<\th_0<0.2$.
Choose a
cut-off function $\psi\in C^\infty_c(\R)$ such that $\psi(\th)=1$ for $|\psi-\th_0| \le \th_0/8$ and $\psi(\th)=0$ for $|\psi-\th_0| > \th_0/4$. Let
$\eta$ be a mollifier supported in $(-1,1)$ and
\[
h(\th) =(1- \psi) g_1 + (\psi g_1)* \eta_{\th_0/8} .
\]
We have $h(\th)$ is smooth for $\th \in (0, \infty)$, and
\[
h(\th)=
g(\th), \quad  (0 \le \th \le \frac 34 \th_0); \qquad
h(\th)=
R_1, \quad  (\frac 54 \th_0  \le \th < \infty).
\]
Modify $h(\th)$ so that it is even with respect to $\pi/2$,
\[
h_1(\th)=
h(\th) ,\quad  (0 \le \th \le \frac \pi 2);\qquad
h_1(\th)=
h(\pi - \th) \quad  (\frac \pi 2  \le \th < \pi) .
\]
The desired domain $V$ is enclosed by the line segment $[-1,1]$ and the curve $r=h_1(\th)$.
\end{proof}
\section*{Acknowledgments}
We warmly thank Hyunwoo Kwon for very
helpful comments and reference \cite{DinMit}. The research of both MC and TT was partially supported by Natural Sciences and Engineering Research Council of Canada (NSERC) under grant RGPIN-2023-04534.

\addcontentsline{toc}{section}{\protect\numberline{}{\hspace{2mm}References}}

\medskip

Misha Chernobai, Department of Mathematics, University of British Columbia, Vancouver, BC V6T 1Z2, Canada;
e-mail: mchernobay@gmail.com

Current affiliation: Max Planck Institute, Leipzig
\medskip

Tai-Peng Tsai, Department of Mathematics, University of British Columbia,
Vancouver, BC V6T 1Z2, Canada; e-mail: ttsai@math.ubc.ca

\end{document}